# On approximation of Markov binomial distributions

AIHUA XIA[1] and MEI ZHANG[2]

[1]*Department of Mathematics and Statistics, The University of Melbourne, Parkville, VIC 3052, Australia. E-mail: [xia@ms.unimelb.edu.au](xia@ms.unimelb.edu.au)*

[2]*School of Mathematical Sciences, Beijing Normal University, Beijing, 100875, P. R. China and Department of Mathematics and Statistics, The University of Melbourne, Parkville, VIC 3052, Australia. E-mail: [meizhang@bnu.edu.cn](meizhang@bnu.edu.cn)*

For a Markov chain $\mathbf{X} = \{X_i, i = 1, 2, \ldots, n\}$ with the state space $\{0, 1\}$, the random variable $S := \sum_{i=1}^n X_i$ is said to follow a Markov binomial distribution. The exact distribution of $S$, denoted $\mathcal{L}S$, is very computationally intensive for large $n$ (see Gabriel [*Biometrika* **46** (1959) 454–460] and Bhat and Lal [*Adv. in Appl. Probab.* **20** (1988) 677–680]) and this paper concerns suitable approximate distributions for $\mathcal{L}S$ when $\mathbf{X}$ is stationary. We conclude that the negative binomial and binomial distributions are appropriate approximations for $\mathcal{L}S$ when $\operatorname{Var} S$ is greater than and less than $\mathbb{E}S$, respectively. Also, due to the unique structure of the distribution, we are able to derive explicit error estimates for these approximations.

## 1. Introduction and the main results

Let $\mathbf{X} = \{X_i, i = 1, 2, \ldots, n\}$ be a Markov chain with the state space $\{0, 1\}$ and transition matrix

$$P = \begin{pmatrix} p_{00} & p_{01} \\ p_{10} & p_{11} \end{pmatrix} = \begin{pmatrix} 1-\alpha & \alpha \\ 1-\beta & \beta \end{pmatrix}, \tag{1.1}$$

where $\alpha, \beta \in (0, 1)$. The distribution of $S := \sum_{i=1}^n X_i$, denoted $\mathcal{L}S$, is well known as the Markov binomial distribution. When $\mathbf{X}$ is stationary and $\alpha = \beta$, $\mathcal{L}S$ degenerates to a binomial distribution. Except for the case $\alpha = \beta$, the exact distribution of $S$ (see Gabriel (1959) and Bhat and Lal (1988)) is very computationally intensive for large $n$ and our interest is in investigating suitable approximate distributions for $\mathcal{L}S$.

It appears that Koopman (1950) and Dobrushin (1961) were among the earliest in the study of limit theory of Markov binomial distributions and the topic was then treated in many articles including Serfling (1975), Wang (1981), Serfozo (1986), He and Xia (1997), Čekanavičius and Mikalauskas (1999), Vellaisamy and Chaudhuri (1999), Barbour and Lindvall (2006), Čekanavičius and Roos (2007). The approximate distributions considered are mostly normal, compound Poisson, translated Poisson or binomial distributions. For instance, when $n\alpha/(1 - \beta + \alpha)$ converges, Wang (1981) proved that for any fixed $k$, $\mathbb{P}(S = k)$ converges to $\mathbb{P}(Y = k)$, where $Y$ is a compound Poisson variable. Barbour and Lindvall (2006) used a translated Poisson distribution to approximate the distribution of a sum of integer-valued random variables whose distributions depend on the state of an underlying Markov chain. Under an aperiodic condition, they established error bounds with respect to the total variation distance, comparable to those found for normal approximation with respect to the weaker Kolmogorov distance. On the other hand, when the first two factorial cumulants of $\mathcal{L}S$ are matched by those of a binomial distribution, Čekanavičius and Roos (2007) demonstrated that the binomial distribution is a suitable approximation for $\mathcal{L}S$ with an approximation error, measured in total variation norm, in the order of $\frac{1}{\sqrt{n}}$. The error estimates in Barbour and Lindvall (2006) and Čekanavičius and Roos (2007) are of the best possible order.

The main purpose of this paper is to find suitable approximate distributions for $\mathcal{L}S$ and provide error bounds as explicit functions of the parameters of the Markov binomial distribution. We will show that the negative binomial and binomial distributions are suitable approximations when $\text{Var}\, S$ is greater than and less than $\mathbb{E}S$, respectively. We employ the celebrated Stein method for binomial (Ehm (1991)) and negative binomial (Brown and Phillips (1999)) approximations and use the unique structure of the Markov binomial distribution to construct a suitable coupling which enables us to specify all of the constants involved in the estimates.

For convenience, from now on, we will assume that $\mathbf{X}$ is stationary. Direct computation ensures that the stationary distribution $\pi$ of $\mathbf{X}$ is

$$p := \pi(1) = \frac{\alpha}{1 - \beta + \alpha}, \qquad \pi(0) = \frac{1 - \beta}{1 - \beta + \alpha}$$

and

$$\mathbb{E}S = np,$$
$$\text{Var}\, S = np(1-p) + nA_0 - A_1 + A_1(\beta - \alpha)^n, \tag{1.2}$$

where

$$A_0 = \frac{2\alpha(1-\beta)(\beta-\alpha)}{(1-\beta+\alpha)^3}, \qquad A_1 = \frac{2\alpha(1-\beta)(\beta-\alpha)}{(1-\beta+\alpha)^4}. \tag{1.3}$$

Note that $\mathbf{X}$ is a stationary positive recurrent Markov chain.

To state the main result, we use $\text{Bi}(m, \theta)$ to stand for the binomial distribution with parameters $m$ and $0 < \theta < 1$. We say that $Y$ follows the negative binomial distribution with parameters $r > 0$ and $0 < q < 1$, denoted by $\text{NB}(r, q)$, if

$$\mathbb{P}(Y = k) = \frac{\Gamma(r+k)}{\Gamma(r)k!} q^r (1-q)^k, \qquad k \in \mathbb{Z}_+ := \{0, 1, 2, \ldots\}.$$

The metric we will use for measuring the approximation errors is the total variation distance defined as

$$d_{\mathrm{TV}}(P,Q) := \sup_{A \subset \mathbb{Z}_+} |P(A) - Q(A)|$$

for probability distributions $P, Q$ on $\mathbb{Z}_+$.

For the Markov chain $\mathbf{X}$ with transition matrix (1.1), we set

$$\mu_1 = \frac{1-\alpha}{\alpha}, \qquad \sigma_1^2 = \frac{1-\alpha}{\alpha^2}, \qquad \mu_2 = \frac{\beta}{1-\beta}, \qquad \sigma_2^2 = \frac{\beta}{(1-\beta)^2},$$

$$C_0 = \frac{|\beta - \alpha|(5 + 43\alpha \vee \beta)}{(1 - \beta \vee \alpha)^2}, \qquad C_1 = \frac{10(\beta \vee \alpha)}{1 - \beta \vee \alpha}, \qquad C_2 = \frac{(1-p)(5 + 23\alpha \vee \beta)}{(1 - \alpha \vee \beta)^2},$$

$$K_1 = \sqrt{5}\sqrt{\frac{\mu_1 + \mu_2 + 2}{\min(1-\alpha, \beta, 1/2)}}, \qquad K_2 = \frac{90(\sigma_1^2 + \sigma_2^2)}{\mu_1 + \mu_2 + 2}.$$

It is worthwhile to note that $\mu_1$ (resp., $\mu_2$) is the mean number of revisits of 0's (resp., 1's) before the Markov chain moves to state 1 (resp., 0), and $\sigma_1^2$ and $\sigma_2^2$ are the variances of the corresponding variables. The main result of the paper is as follows.

**Theorem 1.1.**

1. *If* $\operatorname{Var} S \geq \mathbb{E} S$, *then*

$$d_{\mathrm{TV}}(\mathcal{L}S, \mathrm{NB}(r,q)) \leq C_0 \left[ \frac{2K_1}{\sqrt{n}} + \frac{4K_2}{n} + \beta^{\lfloor n/4 \rfloor} \right], \tag{1.4}$$

*where*

$$r = \frac{(\mathbb{E}S)^2}{\operatorname{Var} S - \mathbb{E} S}, \qquad q = \frac{\mathbb{E} S}{\operatorname{Var} S}$$

*and* $\mathrm{NB}(\infty, 1)$ *is understood as the Poisson distribution with parameter* $\mathbb{E} S$.

2. *If* $\operatorname{Var} S < \mathbb{E} S$, *then*

$$d_{\mathrm{TV}}(\mathcal{L}(S), \mathrm{Bi}(m, \theta))$$
$$\leq \left( \frac{|p-\theta|}{1-\theta} C_1 + \frac{|\beta - \alpha|}{1-\theta} C_2 \right) \left[ \frac{2K_1}{\sqrt{n}} + \frac{4K_2}{n} + (\beta \vee \alpha)^{\lfloor n/4 \rfloor} \right] + \frac{\theta^2(\tilde{m} - m)}{np(1-\theta)}, \tag{1.5}$$

*where*

$$\tilde{m} = \frac{(\mathbb{E}S)^2}{\mathbb{E}S - \operatorname{Var} S}, \qquad m = \lfloor \tilde{m} \rfloor, \qquad \theta = \frac{np}{m}$$

*and* $\lfloor \tilde{m} \rfloor$ *is the integer part of* $\tilde{m}$.

**Remark 1.1.** In practical situations, $\alpha$ and $\beta$ are usually fixed, so the bounds in Theorem 1.1 are of order $\frac{1}{\sqrt{n}}$. The constants $K_i$ and $C_i$ are useful when both $\alpha$ and $\beta$ are a reasonable distance from 0 and 1. If $\alpha$ is close to 0 and $\beta$ is close to 1, then $\mathcal{L}S$ is not unimodal, so one should not expect good approximation by a negative binomial or binomial distribution. On the other hand, when $\alpha$ is close to 1 and $\beta$ is close to 0, $\mathbb{E}S$ is close to $\frac{n}{2}$, but $\operatorname{Var} S$ will be close to 0 for even $n$ and $\frac{1}{4}$ for odd $n$, meaning that we should not expect a good binomial approximation in this case either since the accuracy of approximation is a function of $\operatorname{Var} S$. If both $\alpha$ and $\beta$ are close to 0, then Poisson approximation to $\mathcal{L}S$ (see Barbour *et al.* (1992), Theorem 8.H) is generally sufficient. If both $\alpha$ and $\beta$ are close to 1, one should consider approximating $\mathcal{L}(n-S)$ instead of $\mathcal{L}S$.

***Remark 1.2.*** Except when both $\alpha$ and $\beta$ are very small, Poisson approximation to $\mathcal{L}S$ (see Barbour *et al.* (1992), Theorem 8.H) is inadequate since the error bound of Poisson approximation will not become small when $n$ becomes large.

***Remark 1.3.*** Lemma 2.2, proved in the next section, states that a necessary condition for (1.4) is that $\beta > \alpha$.

***Remark 1.4.*** It is easy to see that if $A_0 > p^2$, then $\operatorname{Var} S > \mathbb{E}S$ for sufficiently large $n$. In this case, as $n \to \infty$, $r \approx \frac{np^2}{A_0 - p^2}$ and $q \approx \frac{p}{p + A_0 - p^2}$.

***Remark 1.5.*** As $n \to \infty$, $m \approx \lfloor \frac{np^2}{p^2 - A_0} \rfloor$ and $\theta \approx p - \frac{A_0}{p} < 1$. Note that if $\alpha = \beta$, $\mathcal{L}S$ degenerates to $\operatorname{Bi}(n, p)$ and $\tilde{m} = m$, so the upper bound of (1.5) becomes 0.

***Remark 1.6.*** Although the estimates in Theorem 1.1 are established for stationary $\mathbf{X}$, since a Markov chain with transition matrix (1.1) and any initial distribution converges exponentially fast to the stationary distribution (see the coupling constructed in the proof of Lemma 2.4), our bounds can be adapted for approximating a Markov binomial distribution with any initial distribution, provided that an error estimate for the difference between the Markov binomial distribution and $\mathcal{L}S$ is added to the upper bounds.

## 2. Preliminary studies of the Markov binomial distribution

To prove Theorem 1.1, we need the following preparation.

**Lemma 2.1.** *Suppose $\{Y_j : j \geq 0\}$ is a Markov chain with transition matrix (1.1) and $Y_0 = 0$. Define $W = \sum_{i=1}^n Y_i$. We then have*

$$d_{\mathrm{TV}}(\mathcal{L}(W), \mathcal{L}(W+1)) \leq \gamma(n), \tag{2.1}$$

*where*

$$\gamma(x) := \frac{K_1}{\sqrt{x}} + \frac{K_2}{x} \qquad \text{for } x > 0,$$

*and $K_1$ and $K_2$ are as given in Section 1.*

**Proof.** We construct another version of the Markov chain $\{Y_i\}$, denoted $\{Y_i'\}$, such that $\mathbb{P}(W + 1 \neq W') \leq \gamma(n)$, where $W' = \sum_{i=1}^n Y_i'$. To this end, let $\rho_0 = 0$ and for $j \geq 1$, let $\rho_j = \inf\{t > \rho_{j-1} : Y_t \neq Y_{\rho_{j-1}}\}$. The $\{\rho_j\}$ are then stopping times separating the Markov chain into blocks of 0's and 1's. In other words, if we set $\xi_j = \rho_j - \rho_{j-1} - 1$ for $j \geq 1$, then $\xi_1$ is the number of revisits of 0's for the Markov chain before it moves to state 1, followed by $\xi_2$ revisits of 1's before it moves to 0, etcetera. By the regenerative theory (see Thorisson (2000), page 53), $\{\xi_j : j \geq 1\}$ are independent random variables, $\xi_{2j-1}$ follows the geometric distribution with parameter $\alpha$ and $\xi_{2j}$ has geometric distribution with parameter $1 - \beta$ for all $j \geq 1$. We write

$$\mu_1 = \mathbb{E}\xi_1 = \frac{1-\alpha}{\alpha}, \qquad \mu_2 = \mathbb{E}\xi_2 = \frac{\beta}{1-\beta},$$

$$\sigma_1^2 = \operatorname{Var}\xi_1 = \frac{1-\alpha}{\alpha^2}, \qquad \sigma_2^2 = \operatorname{Var}\xi_2 = \frac{\beta}{(1-\beta)^2}.$$

For fixed $n$, there are about $\frac{n}{\mu_1+\mu_2+2}$ blocks of 0's and 1's, so we let $k = \lfloor cn \rfloor + 1$ with $c$ close to $(\mu_1 + \mu_2 + 2)^{-1}$. On the other hand, to further simplify the estimate in (2.6) below, it is convenient to take $c = \frac{4}{5}(\mu_1 + \mu_2 + 2)^{-1}$. Let $T_k = \sum_{j=1}^k \xi_{2j-1}$ and $L_k = \sum_{j=1}^k \xi_{2j}$. Using Barbour and Xia (1999), Proposition 4.6, we have

$$d_{\mathrm{TV}}(\mathcal{L}(T_k), \mathcal{L}(T_k + 1)) \leq \frac{1}{\sqrt{cn \min(u_1, 1/2)}},$$

where $u_1 := 1 - d_{\mathrm{TV}}(\mathcal{L}(\xi_1), \mathcal{L}(\xi_1 + 1)) = 1 - \alpha$. We then choose a maximal coupling $(T_k, T_k' + 1)$ of $\mathcal{L}(T_k)$ and $\mathcal{L}(T_k + 1)$ (Barbour *et al.* (1992), page 254) such that

$$d_{\mathrm{TV}}(\mathcal{L}(T_k), \mathcal{L}(T_k + 1)) = \mathbb{P}(T_k \neq T_k' + 1) \leq \frac{1}{\sqrt{cn \min(1-\alpha, 1/2)}} \tag{2.2}$$

and write $\{\xi_{2j-1}', 1 \leq j \leq k\}$ for the i.i.d. random variables satisfying $T_k' = \sum_{j=1}^k \xi_{2j-1}'$. On the other hand, since $\{\xi_{2j}, j \geq 1\}$ play exactly the same role as $\{\xi_{2j-1}, j \geq 1\}$ with 0 and 1 swapped, there exists a maximal coupling $(L_k + 1, L_k')$ of $\mathcal{L}(L_k + 1)$ and $\mathcal{L}(L_k)$ such that $(L_k, L_k')$ is independent of $(T_k, T_k')$ and

$$\mathbb{P}(L_k + 1 \neq L_k') = d_{\mathrm{TV}}(\mathcal{L}(L_k + 1), \mathcal{L}(L_k)) \leq \frac{1}{\sqrt{cn \min(\beta, 1/2)}}. \tag{2.3}$$

We write $\{\xi_{2j}', 1 \leq j \leq k\}$ for the i.i.d. random variables satisfying $L_k' = \sum_{j=1}^k \xi_{2j}'$.

Define $\rho_0' = 0$ and $\rho_j' = \rho_{j-1}' + \xi_j' + 1$, $1 \leq j \leq 2k$. We now couple $\{Y_i'\}$ with $\{Y_i\}$ by setting

$$Y_i' = \begin{cases} 0, & \text{for } \rho_{2j-2}' \leq i < \rho_{2j-1}', 1 \leq j \leq k, \\ 1, & \text{for } \rho_{2j-1}' \leq i < \rho_{2j}', 1 \leq j \leq k, \\ Y_i, & \text{for } i \geq \rho_{2k}'. \end{cases}$$

Under the conditions that $\rho_{2k} \leq n$, $T_k = T'_k + 1$ and $L_k + 1 = L'_k$, we have $W' = W + 1$. Hence,

$$\mathbb{P}(W + 1 \neq W') \leq \mathbb{P}(\rho_{2k} > n) + \mathbb{P}(T_k \neq T'_k + 1) + \mathbb{P}(L_k + 1 \neq L'_k). \quad (2.4)$$

Without loss of generality, we may assume that $cn > 8$. In fact, if $cn \leq 8$, then $\frac{K_1}{\sqrt{n}} \geq 1$ and (2.1) clearly holds. Using Chebyshev's inequality, we get

$$\begin{aligned}
\mathbb{P}(\rho_{2k} > n) &\leq \frac{\operatorname{Var}(\rho_{2k})}{(n - \mathbb{E}\rho_{2k})^2} = \frac{k(\sigma_1^2 + \sigma_2^2)}{(n - k(\mu_1 + \mu_2 + 2))^2} \\
&\leq \frac{(cn+1)(\sigma_1^2 + \sigma_2^2)}{(n - (cn+1)(\mu_1 + \mu_2 + 2))^2} \\
&\leq \frac{1.125cn(\sigma_1^2 + \sigma_2^2)}{(n - 1.125cn(\mu_1 + \mu_2 + 2))^2} \leq \frac{K_2}{n}.
\end{aligned} \quad (2.5)$$

Finally, combining the estimates (2.2), (2.3) and (2.6) with (2.4) yields (2.1). $\square$

**Lemma 2.2.** *If* $\operatorname{Var} S \geq \mathbb{E}S$, *then* $\beta > \alpha$.

**Proof.** By (1.2), we have

$$\begin{aligned}
\operatorname{Var} S - \mathbb{E}S &= -np^2 + nA_0 - A_1(1 - (\beta - \alpha)^n) \\
&= -np^2 + nA_0 - A_0(1 + (\beta - \alpha) + (\beta - \alpha)^2 + \cdots + (\beta - \alpha)^{n-1}) \\
&= -np^2 + A_0(n - (1 + (\beta - \alpha) + (\beta - \alpha)^2 + \cdots + (\beta - \alpha)^{n-1})).
\end{aligned}$$

Clearly, $n - (1 + (\beta - \alpha) + (\beta - \alpha)^2 + \cdots + (\beta - \alpha)^{n-1}) > 0$. If $\beta \leq \alpha$, then $A_0 = \frac{2\alpha(1-\beta)(\beta-\alpha)}{(1-\beta+\alpha)^3} \leq 0$, so $\operatorname{Var} S - \mathbb{E}S < 0$, contradicting the assumption. $\square$

**Lemma 2.3.** *If $h$ is a bounded function on $\mathbb{Z}_+$, and $V_1$, $V_2$ and $V$ are $\mathbb{Z}_+$-valued random variables coupled in such a way that $V$ is independent of $(V_1, V_2)$, then*

$$\left| \mathbb{E} \sum_{j=0}^{V_2 - 1} [h(V_1 + j + V) - h(V)] \right| \leq 2\varepsilon_V \|h\| \mathbb{E}\left[ V_1 V_2 + \frac{1}{2}(V_2 - 1)V_2 \right],$$

*where* $\varepsilon_V := d_{\mathrm{TV}}(\mathcal{L}(V), \mathcal{L}(V + 1))$.

**Proof.** We write $\Delta h(\cdot) = h(\cdot + 1) - h(\cdot)$. Then,

$$\left| \mathbb{E} \sum_{j=0}^{V_2 - 1} [h(V_1 + j + V) - h(V)] \right|$$

$$= \left| \sum_{i_1, i_2} \mathbb{E}\left( \sum_{j=0}^{i_2 - 1} [h(i_1 + j + V) - h(V)] \Big| (V_1, V_2) = (i_1, i_2) \right) \mathbb{P}((V_1, V_2) = (i_1, i_2)) \right|$$

$$\leq \sum_{i_1,i_2} \sum_{j=0}^{i_2-1} |\mathbb{E}[h(i_1+j+V) - h(V)]| \mathbb{P}((V_1,V_2) = (i_1,i_2))$$

$$= \sum_{i_1,i_2} \sum_{j=0}^{i_2-1} \left| \mathbb{E} \sum_{i=0}^{i_1+j-1} \Delta h(V+i) \right| \mathbb{P}((V_1,V_2) = (i_1,i_2))$$

$$\leq 2\varepsilon_V \|h\| \sum_{i_1,i_2} \sum_{j=0}^{i_2-1} \sum_{i=0}^{i_1+j-1} \mathbb{P}((V_1,V_2) = (i_1,i_2))$$

$$= 2\varepsilon_V \|h\| \mathbb{E}\left[V_1 V_2 + \frac{1}{2}(V_2-1)V_2\right]. \qquad \square$$

**Lemma 2.4.** *Write* $\mathcal{L}(S - X_i | X_i = j) := \mathcal{L}(S^{i,j})$ *for* $1 \leq i \leq n$ *and* $j = 0, 1$. *If* $h$ *is a bounded function on* $\mathbb{Z}_+$, *then*

$$|\mathbb{E}h(S^{i,1}) - \mathbb{E}h(S)| \leq \|h\| \frac{10\alpha \vee \beta}{1 - \alpha \vee \beta} (\gamma(n/4) + (\alpha \vee \beta)^{\lfloor n/4 \rfloor}), \qquad (2.6)$$

$$\begin{aligned}
&|\mathbb{E}[h(S^{i,1}) - h(S^{i,0})] - \mathbb{E}(S^{i,1} - S^{i,0})\mathbb{E}\Delta h(S)| \\
&\leq \|\Delta h\| \frac{|\alpha - \beta|(5 + 23\alpha \vee \beta)}{(1 - \alpha \vee \beta)^2} (\gamma(n/4) + (\alpha \vee \beta)^{\lfloor n/4 \rfloor}).
\end{aligned} \qquad (2.7)$$

**Proof.** We construct two copies of Markov chains having transition matrix (1.1), with one starting at state 1 and the other at state 0 at time $i$ in such a way that they can meet as soon as possible in both directions and, once they meet, they stay together from then on. To this end, we define a two-dimensional Markov chain $\{(Z_l^{i,1}, Z_l^{i,0}), l \geq i\}$ with state space $\{(0,0),(0,1),(1,0),(1,1)\}$, initial state $(Z_i^{i,1}, Z_i^{i,0}) = (1,0)$ and transition probabilities

$$\begin{aligned}
p_{(1,0)(j_2,j_1)} = p_{(0,1)(j_1,j_2)} &= \begin{cases} p_{0j} \wedge p_{1j}, & \text{if } j_1 = j_2 = j, \\ \beta - \alpha, & \text{if } \beta > \alpha, j_1 = 0, j_2 = 1, \\ \alpha - \beta, & \text{if } \beta \leq \alpha, j_1 = 1, j_2 = 0; \end{cases} \\
p_{(i,i)(j,j)} &= p_{ij} \qquad \text{for } i,j = 0,1.
\end{aligned} \qquad (2.8)$$

Since the reverse chain $\tilde{\mathbf{X}}$ of $\mathbf{X}$ has the same transition matrix as that of $\mathbf{X}$, we can construct a reverse chain $\{(Z_l^{i,1}, Z_l^{i,0}), l < i\}$ of $\{(Z_l^{i,1}, Z_l^{i,0}), l > i\}$ in the same way as in (2.8).

As $i$ is fixed, we drop the subindex $i$ and define

$$\begin{aligned}
\varsigma &= \min\{t - i > 0 : Z_t^{i,1} = Z_t^{i,0}\}, \\
\tilde{\varsigma} &= \min\{i - t > 0 : Z_t^{i,1} = Z_t^{i,0}\}, \\
\tau &= \min\{t - i > 0 : Z_t^{i,1} = Z_t^{i,0} = 0\}
\end{aligned}$$

and

$$\tilde{\tau} = \min\{i - t > 0 : Z_t^{i,1} = Z_t^{i,0} = 0\}.$$

$\varsigma$ and $\tilde{\varsigma}$ then have the same distribution, as do $\tau$ and $\tilde{\tau}$. Moreover,

$$\mathbb{P}(\varsigma \geq m) = |\beta - \alpha|^{m-1}, \qquad \mathbb{P}(\tau \geq m) \leq (\beta \vee \alpha)^{m-1}, \qquad m \geq 1,$$

$$\mathbb{E}(\tilde{\tau} - 1) = \mathbb{E}(\tau - 1) \leq \frac{\beta \vee \alpha}{1 - \beta \vee \alpha}, \tag{2.9}$$

$$\mathbb{E}[(\tilde{\tau} - 1)^2] = \mathbb{E}[(\tau - 1)^2] \leq \frac{(\beta \vee \alpha)(1 + \beta \vee \alpha)}{(1 - \beta \vee \alpha)^2}.$$

By (2.8) and the regenerative theory, the left range $\{(Z_l^{i,1}, Z_l^{i,0}) : l < i - \tilde{\tau}\}$, the middle range $\{(Z_l^{i,1}, Z_l^{i,0}) : l \in [i - \tilde{\tau}, i + \tau]\}$ and the right range $\{(Z_l^{i,1}, Z_l^{i,0}) : l > i + \tau\}$ are independent. If we stipulate $\sum_b^a = 0$ for $a < b$ and let

$$S_l^i = \sum_{j=1}^{i-\tilde{\tau}-1} Z_j^{i,1}, \qquad S_r^i = \sum_{j=i+\tau+1}^{n} Z_j^{i,1},$$

$$\zeta^{i,1} = \sum_{j=(i-\tilde{\tau})\vee 1, j \neq i}^{(i+\tau)\wedge n} Z_j^{i,1}, \qquad \zeta^{i,0} = \sum_{j=(i-\tilde{\tau})\vee 1, j \neq i}^{(i+\tau)\wedge n} Z_j^{i,0}, \tag{2.10}$$

then we can write

$$S^{i,1} = S_l^i + S_r^i + \zeta^{i,1}, \qquad S^{i,0} = S_l^i + S_r^i + \zeta^{i,0}. \tag{2.11}$$

Let $U_i := S_l^i + S_r^i$. We wish to estimate $\varepsilon_i := d_{\mathrm{TV}}(\mathcal{L}(U_i), \mathcal{L}(U_i + 1))$. Due to the symmetry about $i$ of the Markov chain coupled, it suffices to estimate $\varepsilon_i$ for $i \leq \frac{n}{2}$. By the definition of $S_r^i$ and Lemma 2.1,

$$d_{\mathrm{TV}}(\mathcal{L}(S_r^i + 1), \mathcal{L}(S_r^i))$$

$$= d_{\mathrm{TV}}\left(\mathcal{L}\left(\sum_{j=i+\tau+1}^{n} Z_j^{i,1} + 1\right), \mathcal{L}\left(\sum_{j=i+\tau+1}^{n} Z_j^{i,1}\right)\right)$$

$$\leq \sum_{a \leq n/4} d_{\mathrm{TV}}\left(\mathcal{L}\left(\sum_{j=i+a+1}^{n} Z_j^{i,1} + 1\right), \mathcal{L}\left(\sum_{j=i+a+1}^{n} Z_j^{i,1}\right)\right) \mathbb{P}(\tau = a) + \mathbb{P}(\tau > n/4)$$

$$\leq \gamma(n/4) + \mathbb{P}(\tau > n/4),$$

which, because of the independence of $S_l^i$ and $S_r^i$, ensures that

$$\varepsilon_i \leq \gamma(n/4) + \mathbb{P}(\tau > n/4) \leq \gamma(n/4) + (\alpha \vee \beta)^{\lfloor n/4 \rfloor}. \tag{2.12}$$

To compare $S^{i,1}$, $S^{i,0}$ with $S$, we let $\{Y'_l\} = \{Z^{i,1}_l\}$ with probability $p$ and $\{Y'_l\} = \{Z^{i,0}_l\}$ with probability $1-p$ so that $\{Y'_l : 0 \leq l \leq n\}$ has the same distribution as $\mathbf{X}$. Next, replace $\{Y'_l : l \in [i-\tilde{\tau}, i+\tau]\}$ with $\{Y''_l : l \in [i-\tilde{\tau}, i+\tau]\}$, which has the same distribution as $\{Y'_l : l \in [i-\tilde{\tau}, i+\tau]\}$, but is independent of $\{(Z^{i,1}_l, Z^{i,0}_l) : 1 \leq l \leq n\}$. Define

$$Z'_l = \begin{cases} Y''_l, & l \in [i-\tilde{\tau}, i+\tau], \\ Y'_l, & l > i+\tau \text{ or } l < i-\tilde{\tau}' \end{cases} \quad \text{and} \quad \zeta^i = \sum_{l=(i-\tilde{\tau})\vee 1}^{(i+\tau)\wedge n} Z'_l$$

so that $S' := S^i_l + S^i_r + \zeta^i$ follows the distribution $\mathcal{L}S$. By Lemma 2.3, we have

$$|\mathbb{E}[h(S^{i,1}) - h(S')]|$$
$$\leq |\mathbb{E}[h(U_i + \zeta^{i,1}) - h(U_i)]| + |\mathbb{E}[h(U_i + \zeta^i) - h(U_i)]|$$
$$\leq 2\varepsilon_i \|h\|(\mathbb{E}\zeta^{i,1} + \mathbb{E}\zeta^i).$$

However, it follows from (2.9) that

$$\mathbb{E}(\zeta^{i,1} \vee \zeta^{i,0}) \leq \mathbb{E}(\tau - 1) + \mathbb{E}(\tilde{\tau} - 1) \leq \frac{2\alpha \vee \beta}{1 - \alpha \vee \beta}$$

and

$$\mathbb{E}\zeta^i \leq p\mathbb{E}(\tau + \tilde{\tau} - 1) + (1-p)\mathbb{E}(\tau + \tilde{\tau} - 2) \leq \frac{2\alpha \vee \beta}{1 - \alpha \vee \beta} + p \leq \frac{3\alpha \vee \beta}{1 - \alpha \vee \beta}. \qquad (2.13)$$

Therefore,

$$|\mathbb{E}[h(S^{i,1}) - h(S')]| \leq \frac{10\alpha \vee \beta}{1 - \alpha \vee \beta} \|h\|\varepsilon_i, \qquad (2.14)$$

which, together with (2.12), ensures (2.6).

To estimate (2.7), noting that $\beta > \alpha$ implies $\zeta^{i,1} \geq \zeta^{i,0}$, while $\beta \leq \alpha$ gives $\zeta^{i,1} \leq \zeta^{i,0}$, and swapping 0 and 1 in the superscripts, if necessary, we may assume without loss of generality that $\zeta^{i,1} \geq \zeta^{i,0}$. Observing that $U_i$ is independent of $(\zeta^{i,1} - \zeta^{i,0}, \zeta^i)$, we obtain from Lemma 2.3 that

$$\begin{aligned}
&|\mathbb{E}[h(S^{i,1}) - h(S^{i,0})] - \mathbb{E}(S^{i,1} - S^{i,0})\mathbb{E}\Delta h(S)| \\
&= |\mathbb{E}[h(U_i + \zeta^{i,1}) - h(U_i + \zeta^{i,0})] - \mathbb{E}(\zeta^{i,1} - \zeta^{i,0})\mathbb{E}\Delta h(S')| \\
&\leq |\mathbb{E}[h(U_i + \zeta^{i,1}) - h(U_i + \zeta^{i,0})] - \mathbb{E}(\zeta^{i,1} - \zeta^{i,0})\mathbb{E}\Delta h(U_i)| \\
&\quad + \mathbb{E}(\zeta^{i,1} - \zeta^{i,0})|\mathbb{E}[\Delta h(U_i) - \Delta h(S')]| \\
&\leq \left|\mathbb{E}\sum_{j=0}^{\zeta^{i,1}-\zeta^{i,0}-1}[\Delta h(U_i + \zeta^{i,0} + j) - \Delta h(U_i)]\right| + 2\varepsilon_i \|\Delta h\|\mathbb{E}(\zeta^{i,1} - \zeta^{i,0})\mathbb{E}\zeta^i \\
&\leq 2\varepsilon_i \|\Delta h\|\left\{\mathbb{E}\left[\zeta^{i,0}(\zeta^{i,1} - \zeta^{i,0}) + \frac{1}{2}(\zeta^{i,1} - \zeta^{i,0})(\zeta^{i,1} - \zeta^{i,0} - 1)\right] + \mathbb{E}(\zeta^{i,1} - \zeta^{i,0})\mathbb{E}\zeta^i\right\}.
\end{aligned} \qquad (2.15)$$

Now, again using (2.9), we have

$$\mathbb{E}(\zeta^{i,1} - \zeta^{i,0}) \le \mathbb{E}(\varsigma - 1) + \mathbb{E}(\tilde{\varsigma} - 1) = \frac{2|\alpha - \beta|}{1 - |\alpha - \beta|}, \quad (2.16)$$

$$\mathbb{E}[(\zeta^{i,1} - \zeta^{i,0})(\zeta^{i,1} - \zeta^{i,0} - 1)] \le \mathbb{E}[(\varsigma + \tilde{\varsigma} - 2)(\varsigma + \tilde{\varsigma} - 3)]$$
$$= \frac{6|\alpha - \beta|^2}{(1 - |\alpha - \beta|)^2}. \quad (2.17)$$

To estimate $\mathbb{E}(\zeta^{i,0}(\zeta^{i,1} - \zeta^{i,0}))$, for $\kappa = 0, 1$, define

$$\zeta^{i,\kappa,+} = \sum_{l=i+1}^{(i+\tau)\wedge n} Z_l^{i,\kappa}, \qquad \zeta^{i,\kappa,-} = \sum_{l=(i-\tilde{\tau})\vee 1}^{i-1} Z_l^{i,\kappa}$$

and $\tau_{1,0} = \inf\{l \ge 1 : Z_{i+\varsigma+l}^{i,1} = 0\}$. The conditional distribution of $\tau_{1,0}$ given $Z_{i+\varsigma}^{i,1} = 1$ is then the same as $\mathcal{L}(\xi_2 + 1)$. Since $(\zeta^{i,1,+}, \zeta^{i,0,+})$ and $(\zeta^{i,1,-}, \zeta^{i,0,-})$ are independent, and, for convenience, we may assume that they are identically distributed, it follows that

$$\mathbb{E}[\zeta^{i,0}(\zeta^{i,1} - \zeta^{i,0})] = \mathbb{E}[(\zeta^{i,0,-} + \zeta^{i,0,+})(\zeta^{i,1,+} + \zeta^{i,1,-} - \zeta^{i,0,+} - \zeta^{i,0,-})]$$
$$= 2\mathbb{E}(\zeta^{i,0,+})\mathbb{E}(\zeta^{i,1,-} - \zeta^{i,0,-}) + 2\mathbb{E}[\zeta^{i,0,+}(\zeta^{i,1,+} - \zeta^{i,0,+})]. \quad (2.18)$$

On the other hand,

$$\mathbb{E}(\zeta^{i,0,+}) \le \mathbb{E}(\tau - 1) \le \frac{\beta \vee \alpha}{1 - \beta \vee \alpha}, \quad (2.19)$$

$$\mathbb{E}(\zeta^{i,1,-} - \zeta^{i,0,-}) \le \mathbb{E}(\tilde{\varsigma} - 1) = \frac{|\alpha - \beta|}{1 - |\alpha - \beta|}, \quad (2.20)$$

$$\mathbb{E}[\zeta^{i,0,+}(\zeta^{i,1,+} - \zeta^{i,0,+})] \quad (2.21)$$
$$= \mathbb{E}\left\{\left(\sum_{l=i+1}^{(i+\varsigma-1)\wedge n} Z_l^{i,0} + \sum_{l=(i+\varsigma)\wedge n}^{(i+\tau-1)\wedge n} Z_l^{i,0}\right)\left(\sum_{l=i+1}^{(i+\varsigma-1)\wedge n} (Z_l^{i,1} - Z_l^{i,0})\right)\right\}$$
$$\le \frac{1}{4}\mathbb{E}\left\{\left(\sum_{l=i+1}^{(i+\varsigma-1)\wedge n} Z_l^{i,1}\right)^2\right\}$$
$$+ \sum_{s=0,1} \mathbb{E}\left\{\left(\sum_{l=(i+\varsigma)\wedge n}^{(i+\tau-1)\wedge n} Z_l^{i,0}\right)\left(\sum_{l=i+1}^{(i+\varsigma-1)\wedge n} (Z_l^{i,1} - Z_l^{i,0})\right)\middle| Z_{i+\varsigma}^{i,1} = s\right\}\mathbb{P}(Z_{i+\varsigma}^{i,1} = s)$$
$$\le \frac{1}{4}\mathbb{E}[(\varsigma - 1)^2] + \mathbb{E}[\tau_{1,0}(\varsigma - 1)|Z_{i+\varsigma}^{i,1} = 1]\mathbb{P}(Z_{i+\varsigma}^{i,1} = 1)$$

$$= \frac{1}{4}\mathbb{E}[(\varsigma-1)^2] + \frac{1}{1-\beta}\mathbb{E}[(\varsigma-1)|Z_{i+\varsigma}^{i,1}=1]\mathbb{P}(Z_{i+\varsigma}^{i,1}=1)$$

$$\leq \frac{1}{4}\mathbb{E}[(\varsigma-1)^2] + \frac{1}{1-\beta}\mathbb{E}(\varsigma-1)$$

$$\leq \frac{|\alpha-\beta|(1.25+0.25\alpha\vee\beta)}{(1-\alpha\vee\beta)^2},$$

where $\frac{1}{4}$ in the first inequality of (2.21) is due to the fact that $a(b-a)\leq \frac{b^2}{4}$ for all $a$ and $b$. Combining (2.18)–(2.21) yields

$$\mathbb{E}[\zeta^{i,0}(\zeta^{i,1}-\zeta^{i,0})] \leq \frac{2.5|\alpha-\beta|(1+\alpha\vee\beta)}{(1-\alpha\vee\beta)^2}. \qquad (2.22)$$

Therefore, collecting the estimates of (2.13), (2.16), (2.17) and (2.22), we obtain from (2.15) that

$$|\mathbb{E}[h(S^{i,1})-h(S^{i,0})] - \mathbb{E}(S^{i,1}-S^{i,0})\mathbb{E}\Delta h(S)| \leq \varepsilon_i \|\Delta h\| \frac{|\alpha-\beta|(5+23\alpha\vee\beta)}{(1-\alpha\vee\beta)^2},$$

which, together with (2.12), yields (2.7). $\square$

## 3. Proofs of the main results

**Proof of (1.4).** Set $a = r(1-q)$ and $b = 1-q$. Let

$$\mathcal{B}g(j) = (a+bj)g(j+1) - jg(j)$$

be the Stein operator for the negative binomial distribution $\text{NB}(r,q)$ (Brown and Xia (2001)). For $A \subset \mathbb{Z}_+$, let $g_A: \mathbb{Z}_+ \to \mathbb{R}$ be the bounded solution of the Stein equation

$$\mathcal{B}g(j) = \mathbf{1}_{\{j\in A\}} - \text{NB}(r,q)(A) \qquad \text{for all } j \geq 0.$$

Then,

$$d_{\text{TV}}(\mathcal{L}(S), \text{NB}(r,q)) = \sup_{A\subset\mathbb{Z}_+} |\mathbb{E}\mathbf{1}_{\{j\in A\}}(S) - \text{NB}(r,q)(A)| = \sup_{A\subset\mathbb{Z}_+} |\mathbb{E}\mathcal{B}g_A(S)|.$$

It hence remains to show that $|\mathbb{E}\mathcal{B}g_A(S)|$ is bounded by the right-hand side of (1.4) for every $A \subset \mathbb{Z}_+$. For convenience, we drop the subindex $A$ and write $g$ for $g_A$, and define $g'(\cdot) = g(\cdot+1)$. Brown and Xia (2001), Theorem 2.10, states that

$$\|\Delta g'\| := \sup_{j\in\mathbb{Z}_+} |\Delta g'(j)| \leq \frac{1}{a}. \qquad (3.1)$$

Direct computation gives

$$\mathbb{E}\mathcal{B}g(S) = a\mathbb{E}g(S+1) - (1-b)p\sum_{i=1}^{n}\mathbb{E}g(S^{i,1}+2) + p\sum_{i=1}^{n}\mathbb{E}\Delta g(S^{i,1}+1)$$

$$= a\mathbb{E}g'(S) - (1-b)p\sum_{i=1}^{n}\mathbb{E}g'(S^{i,1}+1) + p\sum_{i=1}^{n}\mathbb{E}\Delta g'(S^{i,1}).$$

Let

$$a = n(1-b)p. \tag{3.2}$$

Then,

$$\mathbb{E}\mathcal{B}g(S) = (1-b)p^2\sum_{i=1}^{n}\mathbb{E}g'(S^{i,1}+1) + (1-b)p(1-p)\sum_{i=1}^{n}\mathbb{E}g'(S^{i,0})$$

$$- (1-b)p\sum_{i=1}^{n}\mathbb{E}g'(S^{i,1}+1) + p\sum_{i=1}^{n}\mathbb{E}\Delta g'(S^{i,1})$$

$$= [p^2 + bp(1-p)]\sum_{i=1}^{n}\mathbb{E}\Delta g'(S^{i,1}) - (1-b)p(1-p)\sum_{i=1}^{n}\mathbb{E}[g'(S^{i,1}) - g'(S^{i,0})].$$

Set

$$n[p^2 + bp(1-p)] = (1-b)p(1-p)\sum_{i=1}^{n}\mathbb{E}(S^{i,1} - S^{i,0}),$$

which is equivalent to

$$1 - b = \frac{\mathbb{E}S}{\operatorname{Var}S}. \tag{3.3}$$

Hence, we can write

$$\mathbb{E}\mathcal{B}g(S) = [p^2 + bp(1-p)]\sum_{i=1}^{n}\mathbb{E}[\Delta g'(S^{i,1}) - \Delta g'(S)] \tag{3.4}$$

$$- (1-b)p(1-p)\sum_{i=1}^{n}\{\mathbb{E}[g'(S^{i,1}) - g'(S^{i,0})] - \mathbb{E}(S^{i,1} - S^{i,0})\mathbb{E}\Delta g'(S)\}.$$

Since $\alpha < \beta$ (see Lemma 2.2), we have

$$\frac{p + b(1-p)}{1-b} = p + \frac{\operatorname{Var}S - \mathbb{E}S}{\mathbb{E}S} \leq \frac{2(\beta - \alpha)}{1 - \alpha \vee \beta},$$

so (1.4) follows from applying Lemma 2.4 and (3.1) in (3.4) and then collecting like terms.

Finally, the constants $a$ and $b$ are determined by (3.2) and (3.3). □

The proof of (1.5) is based on the Stein operator for the binomial distribution $\text{Bi}(m, \theta)$,

$$\mathcal{B}g(j) = \theta(m-j)g(j+1) - (1-\theta)jg(j), \qquad j \in \mathbb{Z}_+$$

(see Ehm (1991) or Barbour *et al.* (1992), page 188). The idea of the proof is similar to that in Soon (1996), but at the cost of a slight increase in complexity, we can achieve the better estimate (1.5). As $\text{Bi}(m, \theta)$ has support on $\{0, 1, \ldots, m\}$ while $S$ has support on $\{0, 1, \ldots, n\}$ and it is possible that $n > m$, in estimating the distance between $\mathcal{L}S$ and $\text{Bi}(m, \theta)$, one often needs to deal with $S$ on $\{S \geq m+1\}$ separately. The following technical lemma helps us to avoid this issue.

**Lemma 3.1.** *For each $A \subset \mathbb{Z}_+$, there exists a bounded function $g_A$ on $\mathbb{Z}_+$ such that*

$$\mathcal{B}g_A(j) \geq \mathbf{1}_{\{j \in A\}} - \text{Bi}(m, \theta)(A) \qquad \text{for } j \in \mathbb{Z}_+ \tag{3.5}$$

*and*

$$\|\Delta g_A\| \leq \frac{1}{m\theta(1-\theta)}. \tag{3.6}$$

**Proof.** For $0 \leq j \leq m$, define $g_A(j)$ as in Barbour *et al.* (1992), page 189, that is, $g_A(j), 0 \leq j \leq m$, is the solution to the Stein equation

$$\mathcal{B}g_A(j) = \mathbf{1}_{\{j \in A\}} - \text{Bi}(m, \theta)(A), \qquad 0 \leq j \leq m. \tag{3.7}$$

For $j \geq m + 1$, let

$$g_A(j) = \begin{cases} -\dfrac{1 - \theta \text{Bi}(m, \theta)(A)}{m\theta(1-\theta)}, & \text{if } m \notin A, \\ -\dfrac{1 + \theta - \theta \text{Bi}(m, \theta)(A)}{m\theta(1-\theta)}, & \text{if } m \in A. \end{cases}$$

Direct verification then ensures that

$$\mathcal{B}g_A(j) \begin{cases} = \mathbf{1}_{\{j \in A\}} - \text{Bi}(m, \theta)(A), & \text{if } 0 \leq j \leq m, \\ \geq 1 - \text{Bi}(m, \theta)(A), & \text{if } j \geq m+1, \end{cases}$$

which, in turn, implies (3.5). Using (3.7) with $j = m$, we conclude that

$$g_A(m) = \begin{cases} \dfrac{\text{Bi}(m, \theta)(A)}{m(1-\theta)}, & \text{if } m \notin A, \\ \dfrac{-1 + \text{Bi}(m, \theta)(A)}{m(1-\theta)}, & \text{if } m \in A. \end{cases}$$

Thus,

$$|\Delta g_A(j)| = \begin{cases} \dfrac{1}{m\theta(1-\theta)}, & \text{if } j = m, \\ 0, & \text{if } j \geq m+1. \end{cases}$$

The claim (3.6) follows easily from the proof of Lemma 9.2.1, Barbour *et al.* (1992). □

**Proof of (1.5).** Let $A_0 := \{i : \mathbb{P}(S = i) \geq \text{Bi}(m, \theta)\{i\}\}$ and abbreviate $g_{A_0}$ to $g$. From Lemma 3.1, we have that

$$d_{\text{TV}}(\mathcal{L}S, \text{Bi}(m, \theta)) = \mathbb{P}(S \in A_0) - \text{Bi}(m, \theta)(A_0) \leq \mathbb{E}\mathcal{B}g(S). \tag{3.8}$$

Therefore, it remains to show that $\mathbb{E}\mathcal{B}g(S)$ is bounded by the right-hand side of (1.5). To this end,

$$\begin{aligned}
\mathbb{E}\mathcal{B}g(S) &= \theta\mathbb{E}[(m - S)g(S+1)] - (1 - \theta)\mathbb{E}[Sg(S)] \\
&= m\theta\mathbb{E}[g(S+1)] - \theta\mathbb{E}[S\Delta g(S)] - \mathbb{E}[Sg(S)] \\
&= p(p - \theta)\sum_{i=1}^{n}\mathbb{E}\Delta g(S^{i,1} + 1) - p(1 - p)\sum_{i=1}^{n}\mathbb{E}[g(S^{i,1} + 1) - g(S^{i,0} + 1)] \\
&= p(p - \theta)\sum_{i=1}^{n}[\mathbb{E}\Delta g(S^{i,1} + 1) - \mathbb{E}\Delta g(S + 1)] \tag{3.9} \\
&\quad - p(1 - p)\sum_{i=1}^{n}\{\mathbb{E}[g(S^{i,1} + 1) - g(S^{i,0} + 1)] - \mathbb{E}(S^{i,1} - S^{i,0})\mathbb{E}\Delta g(S + 1)\} \\
&\quad + \left(np(p - \theta) + p(1 - p)\sum_{i=1}^{n}\mathbb{E}(S^{i,0} - S^{i,1})\right)\mathbb{E}\Delta g(S + 1) \\
&:= I_1 + I_2 + I_3.
\end{aligned}$$

By Lemma 2.4 and (3.6), we have

$$|I_1| \leq 10\frac{|p - \theta|}{1 - \theta} \cdot \frac{\beta \vee \alpha}{1 - \beta \vee \alpha}(\gamma(n/4) + (\alpha \vee \beta)^{\lfloor n/4 \rfloor})$$

and

$$|I_2| \leq \frac{1 - p}{1 - \theta}\frac{|\alpha - \beta|(5 + 23\alpha \vee \beta)}{(1 - \alpha \vee \beta)^2}(\gamma(n/4) + (\alpha \vee \beta)^{\lfloor n/4 \rfloor}),$$

which, in turn, ensure that

$$|I_1| + |I_2| \leq \left(\frac{|p - \theta|}{1 - \theta}C_1 + \frac{|\alpha - \beta|}{1 - \theta}C_2\right)\left(\frac{2K_1}{\sqrt{n}} + \frac{4K_2}{n} + (\beta \vee \alpha)^{\lfloor n/4 \rfloor}\right). \tag{3.10}$$

To estimate $I_3$, setting $\epsilon = \tilde{m} - m$, we get

$$\left| np(p-\theta) + p(1-p)\sum_{i=1}^{n}\mathbb{E}(S^{i,0} - S^{i,1}) \right|$$
$$= |np(p-\theta) + np(1-p) - \operatorname{Var} S|$$
$$= \left| \mathbb{E}S - \operatorname{Var} S - \frac{(\mathbb{E}S)^2}{m} \right| = \frac{\epsilon(\mathbb{E}S - \operatorname{Var} S)}{\tilde{m} - \epsilon}$$
$$= \frac{m}{\tilde{m}}\theta^2 \epsilon.$$

Therefore, recalling $\|\Delta g\| \leq \frac{1}{m\theta(1-\theta)}$ and $\tilde{m} \geq m = \frac{np}{\theta}$, we arrive at

$$|I_3| \leq \frac{m}{\tilde{m}}\|\Delta g\|\theta^2 \epsilon \leq \frac{\theta\epsilon}{\tilde{m}(1-\theta)} \leq \frac{\theta^2 \epsilon}{np(1-\theta)}. \tag{3.11}$$

The proof is completed by combining the estimates (3.8)–(3.11). $\square$

## Acknowledgement

This work was done when the corresponding author Mei Zhang worked as a research fellow at the University of Melbourne, and was supported by the ARC Centre of Excellence for Mathematics and Statistics of Complex Systems.